\documentclass{ws-idaqp}
\usepackage{graphicx}
\usepackage[super]{cite}

\usepackage{amsmath}
\usepackage{amssymb}
\usepackage{amsfonts}
\usepackage{amsfonts}
\usepackage{amssymb}
\usepackage{cases}
\usepackage{bbm}
\usepackage{mathrsfs}
\usepackage{graphicx}
\usepackage{amsfonts}
\usepackage{enumerate}
\usepackage{theorem}
\usepackage[pdfstartview=FitH]{hyperref}

\allowdisplaybreaks
\newtheorem{lem}{Lemma}[section]
\newtheorem{prop}[theorem]{Proposition}
\newtheorem{thm}[theorem]{Theorem}
\newtheorem{cor}[theorem]{Corollary}
\newtheorem{dfn}[theorem]{Definition}

\newtheorem{rem}{Remark}[section]
\def\beq{\begin{equation}}
\def\nneq{\end{equation}}

\def\bthm{\begin{thm}}
\def\nthm{\end{thm}}

\def\blem{\begin{lem}}
\def\nlem{\end{lem}}
\def\bprf{\begin{proof}}
\def\nprf{\end{proof}}
\def\bprop{\begin{prop}}
\def\nprop{\end{prop}}
\def\brmk{\begin{rem}}
\def\nrmk{\end{rem}}

\def\bexa{\begin{exa}}
\def\nexa{\end{exa}}
\def\bcor{\begin{cor}}
\def\ncor{\end{cor}}

\newcommand{\ee}{\mathbb{E}}

\newcommand{\pp}{\mathbb{P}}

\newcommand{\tr}{\mathop{\mathrm{tr}}}

\def\AA{\mathcal A}

\def\FF{\mathcal F}
\def\EE{\mathcal E}

\def\HH{\mathcal H}

\def\e{\varepsilon}

\begin{document}

\markboth{J. Yang and J.L. Zhai}{Asymptotics of stochastic 2D hydrodynamical type systems in unbounded domains}

%%%%%%%%%%%%%%%%%%% Publisher's Area please ignore %%%%%%%%%%%%%%%%%%%%%%
\catchline{}{}{}{}{}
%%%%%%%%%%%%%%%%%%%%%%%%%%%%%%%%%%%%%%%%%%%%%%%%%%%%%%%%%%%%%%%%%%%%%%%%%

\title{ASYMPTOTICS OF STOCHASTIC 2D HYDRODYNAMICAL TYPE SYSTEMS IN UNBOUNDED DOMAINS\\
}

\author{JUAN YANG}

\address{School of Science, Beijing University of posts and telecommunications, \\
No.10 Xitucheng Road, Haidian District, Beijing, 100875, The People's Republican of China\\
yangjuanyj6@gmail.com}

\author{JIANLIANG ZHAI}

\address{School of Mathematical Sciences, University of Science and Technology of China,\\
Key Laboratory of Wu Wen-Tsun Mathematics, Chinese Academy of Sciences,\\
No. 96 Jinzhai Road, Hefei, 230026, The People's Republican of China\\
zhaijl@ustc.edu.cn}

\maketitle

\begin{history}
\received{(Day Month Year)}
\revised{(Day Month Year)}
\published{(Day Month Year)}
\comby{(JIANLIANG ZHAI)}
\end{history}

\begin{abstract}
In this paper, we prove a central limit theorem and establish a moderate deviation principle for 2D stochastic hydrodynamical type systems with multiplicative noise in unbounded domains, which covers 2D Navier-Stokes equations, 2D MHD models and the 2D magnetic $B\acute{e}nard$ problem and also
shell models of turbulence. The weak convergence method plays an important role in obtaining the moderate deviation principle.
\end{abstract}



\section{Introduction}

Let $(\Omega,\mathcal{F},\mathbb{P})$ be a probability space with an increasing family $\{\FF_t\}_{0\le t\le T}$ of the sub-$\sigma$-fields of $\FF$ satisfying the usual conditions.
Consider abstract stochastic evolution
equation of the following form
\begin{equation}\label{hydro01}
\partial_t u + \mathcal{A} u + \mathcal{B}(u,u)+ \mathcal{R}(u)=\sigma(t,u)\dot{W}(t),
\end{equation}
where $W(\cdot)$ is a Wiener process. This abstract nonlinear form covers a host class of 2D hydrodynamical type systems including the following models. (1) When $\mathcal{R}=0$, it is 2D Navier-Stokes equations. (2) 2D magneto-hydrodynamic equations are a combination of the Navier-Stokes equations of fluid dynamics and Maxwell's equations of electromagnetism. Magneto-hydrodynamics is the study of the magnetic properties of electrically conducting fluids like plasmas, liquid metals, salt water or electrolytes. (3) The B$\acute{e}$nard problem is a viscous fluid, in a rectangular container D, which is heated from below and the top surface is taken at constant temperature. 2D Boussinesq model for the B$\acute{e}$nard convection is the Navier-Stokes equations coupled with equations for temperature, i.e. heat equations. (4) 2D magnetic B$\acute{e}$nard problem is the Boussinesq model coupled with magnetic field. (5) 3D Leray $\alpha$-model is quite similar to the Lagrangian averaged Navier-Stokes $\alpha$-model of turbulence. The computational overhead connected with the Leray $\alpha$-model was lower than that of dynamic models, like the Lagrangian averaged Navier-Stokes $\alpha$-model and other sub-grid scale models of turbulence. (6) Shell models of turbulence are useful phenomenological models that retain certain features of the Navier-Stokes equations. These models are the Navier-Stokes equations written in Fourier space containing only local interaction between the modes. Their main computational advantage is the parameterization of the fluctuation of a turbulent field in each octave of wave numbers (called shells) by few representative variables.
We refer to \refcite{CM} and the references therein for examples and more details.
\vskip 0.3cm

Let $u^\e$ denote
the solution of the following equation
\beq\label{hydro09}
du^\e(t)+Au^\e(t)dt+B(u^\e(t))dt+\tilde{R}(t,u^\e(t))dt=\sqrt\e\sigma(t, u^\e(t))dW(t),
\nneq
with the initial condition $u^\e(0)=\xi$ for some fixed point $\xi$ in $H$.
As the parameter $\e$ tends to zero, the solution $u^\e$ of \eqref{hydro09} will tend to the solution $u^0$ of the following deterministic equation
\beq\label{hydro10}
du^0(t)+Au^0(t)dt+B(u^0(t))dt+\tilde{R}(t,u^0(t))dt=0,\ \ \ \text{with}\  u^0(0)=\xi \in H.
\nneq
\vskip 0.1cm
The aim of this paper is to study deviations of $u^\epsilon$ from $u^0$ as $\epsilon\rightarrow 0$. That is, the asymptotic behavior of the trajectory
$$\frac{1}{\sqrt \epsilon \lambda(\e)}(u^\epsilon-u^0)(t),\ \ t\in [0,T].$$
When
$\lambda(\e)=1$, it is related to central limit theorem (CLT for short). And when
$\lambda(\e)$ is the deviation scale verifies
\begin{eqnarray}\label{h}
\lambda(\e)\rightarrow+\infty,\ \ \sqrt \e \lambda(\e)\rightarrow 0, \ as\ \e\rightarrow 0,
\end{eqnarray}
it provides so-called {\it moderate deviation principle} (MDP for short, cf. \refcite{Dembo-Zeitouni}). Throughout the paper, we assume that \eqref{h} is in place.
 \vskip 0.1cm
Moderate deviation is an intermediate estimation between the large deviation with scale $\lambda(\e)=1/\sqrt\e$ and the CLT. Like the large deviations, the moderate deviation theory stems from inferential theory of statistics. In recent years, there is an increasing interest on the study of MDP.
For independent and identically distributed random sequences, Chen \refcite{C1}, \refcite{C2}, \refcite{C3} and Ledoux \refcite{LS} found the necessary and sufficient conditions for MDP. For Markov chain processes, moderate deviation was discussed by Djellout and Guillin \refcite{DG},  Gao \refcite{G1}, \refcite{G2},  and Wu \refcite{W1}, \refcite{W2}.
Baier and Freidlin \refcite{BF} and Guillin \refcite{GA} considered models with averaging.
Moderate deviations for mean field interacting particle models was worked by Douc, Guillin and Najim \refcite{DGN} and Del Moral, Hu and Wu \refcite{DHW}.
Wang and Zhang \refcite{WZ}, Wang, Zhai and Zhang \refcite{WZZ} discussed moderate deviations for stochastic reaction-diffusion equation and 2D stochastic Navier-Stokes equations, respectively. Recently, there are several works on moderate deviations for SPDEs with jump, see \refcite{Budhiraja-Dupuis-Ganguly}, \refcite{DXZZ}, \refcite{WZ2015}.

\vskip 0.3cm

%In recent year, there has been a wide-spread interest in the study of qualitative properties of stochastic models which describe
%cooperative effects in fluids by taking into account macroscopic parameters such as temperature or/and magnetic field.

In this paper, we prove a central limit theorem and establish a moderate deviation principle for the class of abstract nonlinear stochastic models of the form (\ref{hydro09}) with multiplicative noise in unbounded domains, which cover a wide class of mathematical coupled models from fluid dynamics mentioned in first sentence. There are two main difficulties. One is caused by the bilinear term in hydrodynamical systems, we have to deal with the moment of the norm delicately. Another one is that, since we do not assume the compactness of embeddings in the corresponding Gelfand triple $V\subset H\subset V'$, this allows us to cover the important class of hydrodynamical models in unbounded domains. The payoff is that we have to impose some more time regularity assumption on the diffusion coefficient (see (\ref{C6})), much of the problem is caused by this. This is the main
difference from \refcite{WZZ}. The technique introduced by \refcite{CM} will be used in our paper.
%Recently, the second author and his cooperators obtained CLT and MDP for 2D stochastic Navier-Stokes equations in \refcite{WZZ}.

\vskip0.3cm

 The organization of this paper is as follows. In Section 2,
 we shall introduce stochastic 2D hydrodynamical type systems.
We prove the CLT in Section 3.
 Section 4 is devoted to show the MDP.

\vskip0.3cm
Throughout this paper, $C$ is a positive constant
%depending on some parameters $N, f, T, \cdots$,
independent of $\e$ and its value may be different from line to line.

%\section{The Main Text}
%Contributions are to be in English. Authors are encouraged to have
%their contribution checked for grammar. Abbreviations are allowed but
%should be spelt out in full when first used.  Integers ten and below
%are to be spelt out but type as (2+1) dimensions.  Italicize foreign
%language phrases (e.g.~Latin, French).
%
%The text should be in 10 pt Times Roman, single spaced with
%baselineskip of 13~pt. Text area (excluding copyright block and folio)
%is 6.9 inches high and 5 inches wide for the first page.
%Text area (excluding running title) is 7.7 inches high and\break
%5 inches wide for subsequent pages.  Final pagination and
%insertion of running titles will be done by the publisher.
%
\section{Description of the Model}
To formulate the stochastic evolution equation \eqref{hydro09}, we introduce the following standard spaces: let
$(H,|\cdot|,(\cdot,\cdot))$ denote a separable Hilbert space, and $A$ be an self-adjoint
positive linear operator on H.
Set $V=Dom(A^{\frac{1}{2}})$ with norm $\|v\|:=|A^{\frac{1}{2}}v|$, where $v\in V$.
Let $V'$ be the dual of $V$ and $(u,v)$ denote the duality between $u\in V$ and $v \in V'$.  Identifying $H$ with its dual $H'$, we have the dense, continuous embedding $V\hookrightarrow H\cong H'\hookrightarrow V'.$

Consider the following stochastic equation:
\begin{equation}\label{hydro011}
d u(t) + {A} u(t) + B(u(t))+ {\tilde{R}}(t,u(t))=\sigma(t,u(t))dW(t),
\end{equation}
where $B(u):=B(u,u)$.

We assume that the mapping $B:\ V \times V \rightarrow V'$ satisfies the following
antisymmetry and bounded conditions:\\
(C1-1) $B:\ V \times V \rightarrow V'$ is a bilinear continuous mapping.\\
(C1-2) For $u_1, u_2, u_3 \in V$,
\begin{equation}\label{hydro02}
(B(u_1,u_2), u_3) = -(B(u_1,u_3), u_2).
\end{equation}
(C1-3) There exists a Banach interpolation space $\mathcal{H}$
possessing the properties:

\  (i) $V\hookrightarrow  \mathcal{H}  \hookrightarrow V'$;

\  (ii) there exists a constant $a_0>0$ such that for any $v\in V$,
\begin{equation}\label{hydro03}
\|v\|_{\mathcal{H}}^2\leq a_0 |v| \|v\|;
\end{equation}

\ (iii) for every $\eta>0$ there exists $C_{\eta}>0$ such that for
$u_1, u_2, u_3 \in V$,
\begin{equation}\label{hydro04}
|(B(u_1,u_2), u_3) |\leq
\eta \|u_3\|^2 +
C_{\eta} \| u_1  \|_{\mathcal{H}}^2  \| u_2  \|_{\mathcal{H}}^2.
\end{equation}

It is evident that
\eqref{hydro04} is equivalent to the following two inequalities (see Remark 2.1 in \refcite{CM}),
\begin{equation}\label{hydro05}
|(B(u_1,u_2), u_3) |\leq
C_1 \|u_3\|^2 +
C_2 \| u_1  \|_{\mathcal{H}}^2  \| u_2  \|_{\mathcal{H}}^2;
\end{equation}
\begin{equation}\label{hydro06}
|(B(u_1,u_2), u_3) |\leq
C\| u_1  \|_{\mathcal{H}}^2  \|u_2\|   \| u_3  \|_{\mathcal{H}}^2,
\end{equation}
where $C_1,\ C_2,\ C$ are some positive constants.

In view of \eqref{hydro02}, \eqref{hydro03} and \eqref{hydro06}, it yields
\begin{equation}\label{hydro07}
|(B(u_1,u_1), u_2) |\leq
\eta \| u_1  \|^2 +
C_{\eta} | u_1 |^2   \| u_2  \|_{\mathcal{H}}^4.
\end{equation}

And then,
\begin{equation}\label{hydro08}
|(B(u_1)-B(u_2), u_1-u_2) |
=|(B(u_1-u_2), u_2) |
\leq
\eta \| u_1-u_2  \|^2 +
C_{\eta} | u_1-u_2 |^2   \| u_2  \|_{\mathcal{H}}^4.
\end{equation}

The covariance operator $Q$ of the Wiener process $W(\cdot)$ is  a positive symmetric, trace class operator on $H$.
Let $H_0=Q^{1/2}H$. Then $H_0$ is a Hilbert space with the inner product
\beq\label{H0}
\langle u,v\rangle_0=(Q^{-1/2}u,Q^{-1/2}v)\ \ \ \ \forall u,v\in H_0.
\nneq
Let $|\cdot|_0$ denote the norm in $H_0$. Clearly, the embedding of $H_0$ in $H$ is Hilbert-Schmidt, since $Q$ is a trace class operator.
Let $L_Q(H_0;H)$ denote the space of linear operators $S$ such that $SQ^{1/2}$ is a Hilbert-Schmidt operator from $H$ to $H$. Define the norm on the space $L_Q(H_0;H)$ by $|S|_{L_Q}=\sqrt{\tr(SQS^*)}$. Set $L_{(H_0,H)}$ the space of all bounded linear operators from $H_0$ into $H$, and denote
$|\cdot|_{L_{(H_0, H)}}$ its norm.

The noise coefficient $\sigma\in C([0,T]\times V;\ L_Q(H_0;H))$
satisfies the following hypothesis, for all $t\in[0,T],\ u,v\in V$,\\
(C2-1) $|\sigma(t,u)|_{L_{Q}}^2\leq K_0 + K_1 |u|^2$, for some $K_0,\ K_1>0$.\\
(C2-2) $|\sigma(t,u)-\sigma(t,v)|_{L_{Q}}^2 \leq L_1 |u-v|^2,$ for some $L_1>0$.

We introduce another coefficient $\tilde{R}\in C([0,T]\times H;\ H)$ which satisfy, for some $R_0,\ R_1>0$,\\
%C3-1) $|\tilde{\sigma}(t,u)|_{L_{H_0, H}}^2\leq \tilde{K}_0 + \tilde{K}_1 |u|^2
%+ \tilde{K}_{\mathcal{H}} \| u  \|_{\mathcal{H}}^2    $, for some $\tilde{K}_0,\ \tilde{K}_1,\ \tilde{K}_{\mathcal{H}} >0$.
%
%\ \ \ \ \ $|\tilde{\sigma}(t,u)  -  \tilde{\sigma}(t,v) |_{L_{(H_0, H)}}^2\leq
% \tilde{L}_1 |u-v|^2
%+ \tilde{L}_2 \| u -v \|^2    $, for some $\tilde{L}_1,\ \tilde{L}_2 >0$,\\
%where $|\cdot|_{L_{H_0, H}}$ denotes the norm in the space of ${L_{H_0, H}}$ of all bounded linear operators from $H_0$ into $H$.\\
(C3) $|\tilde{R}(t,0)| \leq R_0,$\ \  $|\tilde{R}(t,u) -  \tilde{R}(t,v)| \leq R_1 |u-v|,$ $\forall t\in[0,T],\ u,v\in H$.
\vskip0.3cm

From the Theorem 2.4 in \refcite{CM}, the following lemma holds.
\blem\label{Lem 1}{\rm Assume that (C1-1)-(C1-3), (C2-1)-(C2-2) and (C3) hold.
There exists a constant $\e_0>0$ such that for any $0<\e\le \e_0$, \eqref{hydro09} has a unique weak solution $u^\e$ in $L^2(\Omega;C([0,T];H))\cap L^2(\Omega\times[0,T];V)$. Furthermore, for some constant $C$ independent of $\e$,
\beq
\ee\left(\sup_{0\le t \le T}|u^\e(t)|^4+ \int_0^T\|u^\e(s)\|^2ds
+\int_0^T|u^\e(s)|^2\cdot\|u^\e(s)\|^2ds
+\int_0^T\|u^\e(s)\|_{\mathcal{H}}^4ds
\right)\le C(1+|\xi|^4);
\nneq
and particulary,
\beq
\sup_{0\le t \le T}|u^0(t)|^4+ \int_0^T\|u^0(s)\|^2ds
+\int_0^T|u^0(s)|^2\cdot\|u^0(s)\|^2ds
+\int_0^T\|u^0(s)\|_{\mathcal{H}}^4ds
\le C(1+|\xi|^4),
\nneq
}
here $u^0$ is the solution of (\ref{hydro10}).
\nlem

\section{Central Limit Theorem}
In this section, we will establish the CLT.
\vskip0.3cm

The following result is concerned with the convergence of $u^\e$ as $\e\to 0$, and can be obtained similarly as Proposition 3.1 in \refcite{WZZ}. The proof is omitted.
\bprop \label{u convergence} {\rm Under the conditions (C1-1)-(C1-3), (C2-1)-(C2-2) and (C3), there exists a constant $\e_0>0$ such that, for all $0<\e\le\e_0$,
\beq
\ee\left(\sup_{0\le t \le T}|u^\e(t)-u^0(t)|^2+ \int_0^T\|u^\e(s)-u^0(s)\|^2ds\right)\le
 \e c_{T, K_0,K_1}.
\nneq
}
\nprop

\vskip0.3cm

Let $V^0$ be the solution of the following SPDE:
\beq\label{V0}
dV^0(t)+AV^0(t)dt+B(V^0(t), u^0(t))+B(u^0(t), V^0(t))dt+\tilde{R}'(t,u^0(t)) V^0(t)dt=\sigma(t, u^0(t))dW(t),
\nneq
with initial value $V^0(0)=0$, where $\tilde{R}'$ given in (C4).
For the existence and uniqueness of the solution for \eqref{V0}, we need the following additional assumption
 about $\tilde{R}'$.
 \vskip 0.2cm

\textbf{(C4)} $\tilde{R}':\ [0,T]\times H\rightarrow L(H)$ is Fr\'{e}chet derivative of $\tilde{R}$
w.r.t. the second variable and $\tilde{R}'$ is continuous such that  $|\tilde{R}'(t,u)|_{L(H)}\le  \tilde{R}_0'|u|+\tilde{R}_1'$ for some positive  $\tilde{R}_0'$ and $\tilde{R}_1'$.

\vskip 0.2cm
Similarly as the proof in Theorem 2.4 in \refcite{CM},
the following lemma holds.
\blem\label{Lem 2}{\rm Assume that (C1-1)-(C1-3), (C2-1)-(C2-2) and (C4) hold.
\eqref{V0} has a unique weak solution $V^0$ in $L^2(\Omega;C([0,T];H))\cap L^2(\Omega\times[0,T];V)$. Furthermore, the solution has the following estimate, for some constant $C>0$,
\beq
\ee\Big(\sup_{0\le t\le T}|V^0(t)|^4+\int_0^T\|V^0(s)\|^2ds+\int_0^T|V^0(s)|^2\|V^0(s)\|^2ds+\int_0^T\|V^0(s)\|_{\mathcal{H}}^4ds\Big)\le C(1+|\xi|^4);
\nneq
}
\nlem

\vskip0.3cm
Our first main result is the CLT. We add the following assumption on
$\tilde{R}'$.
\vskip 0.2cm

\textbf{(C5)} There exists $C>0$ such that, for every $u_1,\ u_2\in H$,
 $$|\tilde{R}'(t,u_1)-\tilde{R}'(t,u_2)|_{L(H)}\le C|u_1-u_2|.$$

\bthm[Central Limit Theorem]\label{V con.}
{\rm
Under the conditions (C1-1)-(C1-3), (C2-1)-(C2-2) and (C3)-(C5),  $(u^\e-u^0)/\sqrt\e$ converges to $V^0$ in the space $C([0,T];H)\cap L^2([0,T];V)$ in probability, that is, as $\e\rightarrow0$,
\begin{align}\label{eq clt}
\sup_{0\le t\le T}\left|\frac{u^\e(t )-u^0(t)}{\sqrt\e}-V^0(t)\right|^2+ \int_0^{T}
\left\|\frac{u^\e(s )-u^0(s)}{\sqrt\e}-V^0(s)\right\|^2ds\rightarrow 0,\ \text{in\ probability.}
\end{align}
}
\nthm

\begin{rem}
Because of the property of $\tilde{R}'$, we can only obtain (\ref{eq clt}), which is weaker than Theorem 3.2 in \refcite{WZZ}, i.e.
\begin{align*}
\lim_{\e\rightarrow0}\mathbb{E}\Big(\sup_{0\le t\le T}\left|\frac{u^\e(t )-u^0(t)}{\sqrt\e}-V^0(t)\right|^2+ \int_0^{T}
\left\|\frac{u^\e(s )-u^0(s)}{\sqrt\e}-V^0(s)\right\|^2ds\Big)=0.
\end{align*}
\end{rem}

\bprf\ \
Let $
V^\e(t):=(u^\e(t)-u^0(t))/\sqrt\e.
$ Then $V^\e(0)=0$,
\begin{align}\label{hydro112}
&dV^\e(t)+ AV^\e(t)dt+B(V^\e(t),u^\e(t))dt+B(u^0(t),V^\e(t))dt\notag\\
&+ \frac{1}{\sqrt{\e}}[\tilde{R}(t,u^\e(t))-\tilde{R}(t,u^0(t))]
=\sigma(t,u^\e(t))dW(t),
\end{align}
and
\begin{align*}
&d(V^\e(t)-V^0(t))+  A(V^\e(t)-V^0(t))dt\\&+B(V^\e(t)-V^0(t),u^0(t))dt+B(V^\e(t), u^\e(t)-u^0(t))dt+B(u^0(t),V^\e(t)-V^0(t))dt\\
&+ \frac{1}{\sqrt{\e}}[\tilde{R}(t,u^\e(t))-\tilde{R}(t,u^0(t))]dt-\tilde{R}'(t,u^0(t)) V^0(t)dt\\
=&[\sigma(t,u^\e(t))-\sigma(t,u^0(t))]dW(t).
\end{align*}
By It\^{o}'s formula for $|V^\e(t)-V^0(t)|^2$ and \eqref{hydro02}, we have
\begin{align*}
&d|V^\e(t)-V^0(t)|^2+2\|V^\e(t)-V^0(t)\|^2dt\\
=&-2\Big(B(V^\e(t)-V^0(t),u^0(t)), V^\e(t)-V^0(t)\Big)dt\\
&-2\Big(B(V^\e(t), u^\e(t)-u^0(t)), V^\e(t)-V^0(t)\Big)dt\\
&-2\Big( \frac{1}{\sqrt{\e}}[\tilde{R}(t,u^\e(t))-\tilde{R}(t,u^0(t))]-\tilde{R}'(t,u^0(t)) V^0(t), V^\e(t)-V^0(t)\Big)dt\\
 &+ 2\Big([\sigma(t,u^\e(t))-\sigma(t,u^0(t))]dW(t), V^\e(t)-V^0(t)\Big)\\
 &+|\sigma(t,u^\e(t))-\sigma(t,u^0(t))|^2_{L_Q}dt.
 \end{align*}
 Defining $\tau_N=\inf\{t:|V^\e(t)-V^0(t)|^2+\int_0^t\|V^\e(s)-V^0(s)\|^2ds>N\}$, we have
\begin{align}\label{II0}
&|V^\e(t\wedge\tau_N)-V^0(t\wedge\tau_N)|^2+2 \int_0^{t\wedge\tau_N}\|V^\e(s)-V^0(s)\|^2ds\notag\\
\le& 2\int_0^{t\wedge\tau_N}\left|\Big(B(V^\e(s)-V^0(s), u^0(s)), V^\e(s)-V^0(s)\Big)\right|ds\notag\\
&+2\int_0^{t\wedge\tau_N}\left|\Big(B(V^\e(s),u^\e(s)-u^0(s)), V^\e(s)-V^0(s)\Big)\right|ds\notag\\
&+2\int_0^{t\wedge\tau_N}\left|\Big(\frac{1}{\sqrt{\e}}[\tilde{R}(s,u^\e(s))-\tilde{R}(s,u^0(s))]
-\tilde{R}'(s,u^0(s)) V^0(s), V^\e(s)-V^0(s)\Big)\right|ds\notag\\
&+2\left|\int_0^{t\wedge\tau_N}\Big([\sigma(s,u^\e(s))-\sigma(s,u^0(s))]dW(s), V^\e(s)-V^0(s)\Big)\right|\notag\\
&+\int_0^{t\wedge\tau_N}|\sigma(s,u^\e(s))-\sigma(s,u^0(s))|^2_{L_Q}ds\notag\\
 =:&\sum\limits^{5}_{k=1}I_k.
\end{align}

Because of \eqref{hydro02} and \eqref{hydro07}, we have
\begin{align}\label{II1}
I_1(t)
&\le 2\eta \int_0^{t\wedge\tau_N}\|V^{\e}(s)-V^0(s)\|^2ds
+2C_{\eta}\int_0^{t\wedge \tau_N}|V^{\e}(s)-V^0(s)|^2\cdot|u^0(s)|_{\mathcal{H}}^4ds.
\end{align}

Since $V^{\e}(s)=(u^\e(s)-u^0(s))/\sqrt\e$, \eqref{hydro03} and \eqref{hydro04} yield
\begin{align}\label{II2}
I_2(t)
&\le 2\sqrt\e   \int_0^{t\wedge\tau_N}\eta  \|V^{\e}(s)-V^0(s)\|^2ds
+2 \sqrt\e   C_{\eta}\int_0^{t\wedge \tau_N}|V^{\e}(s)|_{\mathcal{H}}^4ds\notag\\
&\le 2\sqrt\e \eta   \int_0^{t\wedge\tau_N} \|V^{\e}(s)-V^0(s)\|^2ds
+2 \sqrt\e   C_{\eta} a_0 \int_0^{t\wedge \tau_N}|V^{\e}(s)|^2\|V^{\e}(s)\|^2ds.
\end{align}

In view of (C4), (C5) and
\begin{align*}
&\frac{1}{\sqrt\e} [\tilde{R}(s,u^\e(s))-\tilde{R}(s,u^0(s))]
-\tilde{R}'(s,u^0(s))  V^0(s)\\
&=\frac{1}{\sqrt\e} [\tilde{R}(s,u^\e(s))-\tilde{R}(s,u^0(s))]
+\tilde{R}'(s,u^0(s))
\left( [ V^{\e}(s)-V^0(s)]-  \frac{1}{\sqrt\e} [u^{\e}(s)-u^{0}(s)]\right)\\
&=\tilde{R}'(s,u^0(s)) [ V^{\e}(s)-V^0(s)]+\frac{1}{\sqrt\e} \left(\tilde{R}(s,u^\e(s))-\tilde{R}(s,u^0(s))
-\tilde{R}'(s,u^0(s)) [u^{\e}(s)-u^{0}(s)]\right)\\
&=
\tilde{R}'(s,u^0(s)) [ V^{\e}(s)-V^0(s)]\\
&+
\frac{1}{\sqrt\e}\int_0^1 \Big[\Big(\tilde{R}'\Big(s,x(u^\e(s)-u^0(s))+u^0(s)\Big)-\tilde{R}'(s,u^0(s))\Big) [u^{\e}(s)-u^{0}(s)]\Big]dx,
\end{align*}
we obtain
\begin{align}\label{II22}
I_3(t)
&\le 2\int_0^{t\wedge\tau_N}(\tilde{R}_1'+\tilde{R}_0'|u^0(s)|) | V^{\e}(s)-V^0(s)|^2 ds
+
\frac{C}{\sqrt\e}\int_0^{t\wedge\tau_N}| V^{\e}(s)-V^0(s)||u^{\e}(s)-u^0(s)|^2ds\notag\\
&\le 2\int_0^{t\wedge\tau_N}(\tilde{R}_1'+\tilde{R}_0'|u^0(s)|) | V^{\e}(s)-V^0(s)|^2 ds
+\int_0^{t\wedge\tau_N}| V^{\e}(s)-V^0(s)|^2 ds\notag\\
&+\frac{C^2}{\e}\int_0^{t\wedge\tau_N}|u^{\e}(s)-u^0(s)|^4ds.
\end{align}

For any $\delta>0$, introduce sets $A_{\delta}^t=\{ \omega:\ \sup_{0\le s\le t}  |u^{\e}(s)-u^0(s)|\le \delta \}$. It is evident that  $A_{\delta}^t$ is non-increasing w.r.t. $t$.
Burkholder-Davis-Gundy inequality and (C2-2)
imply
\begin{align}\label{II3}
&\ee \Big( \sup_{0\le s\le t}[\mathbf{1}_{A_{\delta}^s }I_4(s)]\Big)\notag\\
&\le 4\ee\left( \int_0^{t\wedge\tau_N}\mathbf{1}_{A_{\delta}^s }|\sigma(s,u^\e(s))
-\sigma(s,u^0(s))|_{L_Q}^2\cdot |V^\e(s)-V^0(s)|^2ds  \right)^{\frac12}\notag\\
&\le 4 \ee\left(\sup_{0\le s\le t\wedge\tau_N}\Big[\mathbf{1}_{A_{\delta}^s }|V^\e(s)-V^0(s)|\Big]\cdot\Big(\int_0^{t\wedge \tau_N}L_1 |u^\e(s)-u^0(s)|^2  ds  \Big)^{\frac12}\right)\notag\\
&\le \frac12\ee\left(\sup_{0\le s\le t}\Big[\mathbf{1}_{A_{\delta}^s }|V^\e(s\wedge\tau_N)-V^0(s\wedge\tau_N)|^2\Big]\right)
+8L_1\ee\Big(\int_0^{t\wedge\tau_N}|u^\e(s)-u^0(s)|^2ds\Big).
\end{align}

From (C2-2), we get
\begin{align}\label{II4}
 I_5(t)
&\le L_1 \int_0^{t\wedge\tau_N}|u^\e(s)-u^0(s)|^2ds.
\end{align}

Multiply by $\mathbf{1}_{A_{\delta}^s }$, taking the supremum up to time $t$ in \eqref{II0}, and then taking the expectation, one obtains
\begin{align}\label{II}
&\frac12 \ee\left(\sup_{0\le s\le t}\Big[\mathbf{1}_{A_{\delta}^s }|V^\e(s\wedge\tau_N)-V^0(s\wedge\tau_N)|^2\Big]\right)\notag\\
&
+ 2(1-\eta-\sqrt\e \eta)\ee\left(\sup_{0\le s\le t}
\Big[\mathbf{1}_{A_{\delta}^s }\int_0^{s\wedge\tau_N}\|V^\e(l)-V^0(l)\|^2dl\Big]\right)
\notag\\
\le& 2C_{\eta}
\ee \left(
    \int_0^{t}\sup_{0\le l\le s }\Big(\mathbf{1}_{A_{\delta}^l }|V^{\e}(l\wedge\tau_N)-V^0(l\wedge\tau_N)|^2\Big)\cdot|u^0(s)|_{\mathcal{H}}^4ds
\right)\notag\\
&+ 2 \sqrt\e   C_{\eta} a_0 \ee \left(\sup_{0\le s\le t}\Big(\mathbf{1}_{A_{\delta}^s } \int_0^{s\wedge \tau_N}|V^{\e}(l)|^2\|V^{\e}(l)\|^2dl\Big)
\right)
\notag\\
& +\ee \left( \int_0^{t\wedge\tau_N} \mathbf{1}_{A_{\delta}^s } \|V^{\e}(s)-V^0(s)\|^2ds \right)
 +2\ee \left( \int_0^{t\wedge\tau_N}  \mathbf{1}_{A_{\delta}^s }(\tilde{R}_1'+\tilde{R}_0'|u^0(s)|) | V^{\e}(s)-V^0(s)|^2 ds\right)\notag\\
&+\frac{C^2}{\e}\ee \left(\sup_{0\le s\le t} \Big(\mathbf{1}_{A_{\delta}^s }\int_0^{s\wedge\tau_N} |u^{\e}(l)-u^{0}(l)|^4dl\Big)
\right) +8L_1\ee\Big(\int_0^{t\wedge\tau_N}|u^\e(s)-u^0(s)|^2ds\Big)\notag\\
& + L_1\ee \left( \int_0^{t\wedge\tau_N}|u^\e(s)-u^0(s)|^2ds
\right)\notag\\
\le& \int_0^{t}  \ee \left(\sup_{0\le l\le s}\Big[\mathbf{1}_{A_{\delta}^l }|V^{\e}(l\wedge\tau_N)-V^0(l\wedge\tau_N)|^2\right)
\cdot\Big[2C_{\eta} |u^0(s)|_{\mathcal{H}}^4+1+2 (\tilde{R}_1'+\tilde{R}_0'|u^0(s)|) \Big]ds
\notag\\
&+ 2 \sqrt\e   C_{\eta} a_0 \ee \left(\sup_{0\le s\le t}\Big(\mathbf{1}_{A_{\delta}^s } \int_0^{s\wedge \tau_N}|V^{\e}(l)|^2\|V^{\e}(l)\|^2dl\Big)
\right)
\notag\\
&+C^2\delta^2\ee \left(\int_0^{t}\frac{1}{\e}  |u^{\e}(s)-u^{0}(s)|^2ds
\right) + 9L_1\ee \left( \int_0^{t}|u^\e(s)-u^0(s)|^2ds
\right)\notag\\
\le& \int_0^{t} \ee \left(\sup_{0\le l\le s}\mathbf{1}_{A_{\delta}^l }|V^{\e}(l\wedge\tau_N)-V^0(l\wedge\tau_N)|^2\right)
\cdot\Big[2C_{\eta} |u^0(s)|_{\mathcal{H}}^4+1+2(\tilde{R}_1'+\tilde{R}_0'|u^0(s)|) \Big]ds
\notag\\
&+ 2 \sqrt\e   C_{\eta} a_0\cdot C
+C^2\delta^2
 + C\e,
\end{align}
in view of
Proposition \ref{u convergence} and Lemma \ref{V4} below. Choosing $\eta+\sqrt\e \eta<1/2$,
applying the Gronwall's inequality to $\ee\left(\sup_{0\le s\le t}\Big[\mathbf{1}_{A_{\delta}^s }|V^\e(s\wedge\tau_N )-V^0(s\wedge\tau_N)|^2\Big]\right)$ and taking $N\rightarrow\infty$, it yields
\begin{align}
&\ee\left(\sup_{0\le s\le t}\Big[\mathbf{1}_{A_{\delta}^s }|V^\e(s )-V^0(s)|^2\Big]\right)
+ \ee\left(\sup_{0\le s\le t}
\Big[\mathbf{1}_{A_{\delta}^s }\int_0^{s}\|V^\e(l)-V^0(l)\|^2ds\Big]\right)
\notag\\
 \le&   C(\sqrt\e
+\delta^2
 + \e)\exp\left(\int_0^T\Big[2C_{\eta} |u^0(s)|_{\mathcal{H}}^4+1+2(\tilde{R}_1'+\tilde{R}_0'|u^0(s)|) \Big]ds\right).
\end{align}
Thus,
for any $l>0$,
\begin{align*}
&\pp \big(\sup_{0\le t\le T}|V^\e(t)-V^0(t)|^2+ \int_0^{T}
\|V^\e(s )-V^0(s)\|^2ds>l\big)
\notag\\
 \le&   \frac1l
 \ee\left(\sup_{0\le s\le T}\Big[\mathbf{1}_{A_{\delta}^s }|V^\e(s )-V^0(s)|^2\Big]\right)\notag\\
&
+ \frac1l \ee\left(\sup_{0\le s\le T}
\Big[\mathbf{1}_{A_{\delta}^s }\int_0^{s}\|V^\e(l)-V^0(l)\|^2ds\Big]\right)
+\pp\big((A_{\delta}^T)^c\big)
\notag\\
\le &  \frac1l C(\sqrt\e
+\delta^2
 + \e)+\pp\big((A_{\delta}^T)^c\big),
\end{align*}
where $(A_{\delta}^T)^c$ denotes the complement of $A_{\delta}^T$ and we have used the fact that $A_\delta^t$ is non-increasing
\emph{w.r.t.} $t$.
By the definition of $A_{\delta}^t$ and Proposition \ref{u convergence}, $\pp\big((A_{\delta}^T)^c\big)
\rightarrow 0,$ as $\e\rightarrow0.$
The arbitrary of $\delta$ implies that $\pp \big(\sup_{0\le t\le T}|V^\e(t)-V^0(t)|^2+ \int_0^{T}
\|V^\e(s )-V^0(s)\|^2ds>l\big)\rightarrow 0,$ as $\e\rightarrow0$
which completes the proof.
\nprf
\vskip0.3cm
For the integrity in the proof of last theorem, we have to show the following estimate. Because the operator $\tilde{R}'$ satisfies
(C4) and (C5), using the similar arguments as the proof of Lemma 3.2 in \refcite{WZZ}, we have
\blem\label{V4}{\rm Under the conditions (C1-1)-(C1-3), (C2-1)-(C2-2) and (C3)-(C5), there exists a constant $\e_0>0$ such that
\beq
\sup_{0< \e\le\e_0}\ee\Big(\int_0^T |V^\e(s)|^2\cdot\|V^\e(s)\|^2ds\Big)<\infty.
\nneq
}
\nlem

\section{Moderate deviations}

 Let $Z^\e=(u^\e-u^0)/(\sqrt\e \lambda(\e))$. Then $Z^\e$ satisfies the following SPDE:
\begin{align}\label{eq main01}
&dZ^{\e}(t)+AZ^{\e}(t)dt+B\Big(Z^{\e}(t), u^0(t)+\sqrt \e \lambda(\e)Z^\e(t)\Big)dt+B(u^0(t), Z^\e(t)) dt\notag\\
&+(\sqrt{\e} \lambda(\e))^{-1}\Big[\tilde{R}\Big(t,u^0(t)+\sqrt \e \lambda(\e)Z^\e(t)\Big)-\tilde{R}(t,u^0(t))\Big]dt\notag\\
&=\lambda^{-1}(\e)\sigma\Big(t, u^0(t)+\sqrt\e \lambda(\e)Z^{\e}(t)\Big)dW(t),
\end{align}
with initial value $Z^\e(0)=0$.  This equation admits a unique solution $Z^\e=\Gamma^{\e}(W(\cdot))$, where $\Gamma^\e$ stands for the solution
functional from $C([0,T];H)$ into $C([0,T];H)\cap L^2([0,T];V)$.

 In this part, we will prove that $Z^\e$ satisfies an LDP on $C([0,T];H)\cap L^2([0,T];V)$ with $\lambda(\e)$ satisfying \eqref{h}.
This special type of LDP is usually called the MDP of $u^\e$ (cf. \refcite{Dembo-Zeitouni}).

\vskip0.3cm
Firstly we recall the general criteria for a large deviation principle (LDP) given in \refcite{Budhiraja-Dupuis}.
Let $\mathcal{E}$ be a Polish space with the Borel $\sigma$-field $\mathcal{B}(\mathcal{E})$.
    \begin{dfn}\label{Dfn-Rate function}
       \emph{\textbf{(Rate function)}} A function $I: \mathcal{E}\rightarrow[0,\infty]$ is called a rate function on
       $\mathcal{E}$,
       if for each $M<\infty$, the level set $\{x\in\mathcal{E}:I(x)\leq M\}$ is a compact subset of $\mathcal{E}$.
         \end{dfn}
    \begin{dfn}
       \emph{\textbf{(LDP)}} Let $I$ be a rate function on $\mathcal{E}$.  A family
       $\{X^\e\}$ of $\EE$-valued random elements is said to satisfy the LDP on $\mathcal{E}$
       with rate function $I$, if the following two conditions
       hold.
       \begin{itemize}
         \item[$(a)$](Large deviation upper bound) For each closed subset $F$ of $\mathcal{E}$,
              $$
                \limsup_{\e\rightarrow 0}\e\log\mathbb{P}(X^\e\in F)\leq- \inf_{x\in F}I(x).
              $$
         \item[$(b)$](Large deviation lower bound) For each open subset $G$ of $\mathcal{E}$,
              $$
                \liminf_{\e\rightarrow 0}\e\log\mathbb{P}(X^\e\in G)\geq- \inf_{x\in G}I(x).
              $$
       \end{itemize}
    \end{dfn}

\vskip0.3cm
Next, we introduce the Skeleton Equations.
The Cameron-Martin space associated with the Wiener process $\{W(t), t\in[0,T]\}$ is given by
\beq\label{Cameron-Martin}
\HH_0:=\left\{h:[0,T]\rightarrow H_0; h \  \text{is absolutely continuous and } \int_0^T|\dot h(s)|_0^2ds<+\infty\right\}.
\nneq
The space $\HH_0$ is a Hilbert space with inner product
 $
 \langle h_1, h_2\rangle_{\HH_0}:=\int_0^T\langle \dot h_1(s), \dot h_2(s)\rangle_0ds.
 $
Let $\AA$ denote  the class of $H_0$-valued $\{\FF_t\}$-predictable processes $\phi$ belonging to $\HH_0$ a.s..
Let $S_N=\{h\in \HH_0; \int_0^T|\dot h(s)|_0^2ds\le N\}$. The set $S_N$ endowed with the weak topology is a Polish space.
Define $\AA_N=\{\phi\in \AA;\phi(\omega)\in S_N, \mathbb{P}\text{-a.s.}\}$.

For any $h\in \HH_0$, consider the deterministic integral equation
\begin{align}\label{eq skeleton}
&dX^h(t)+\Big(AX^h(t)+B(X^h(t),u^0(t))+B(u^0(t),X^h(t))\Big)dt\notag\\
&+\tilde{R}'(t,u^0(t))\cdot  X^h(t)dt
=\sigma(t, u^0(t))\dot h(t)dt,
\end{align}
with initial value $X^h(0)=0$
and for any $\phi^\e\in \AA$, consider
\begin{align}\label{eq main}
&dX^{\e}(t)+AX^{\e}(t)dt+B\Big(X^{\e}(t), u^0(t)+\sqrt \e \lambda(\e)X^\e(t)\Big)dt+B(u^0(t), X^\e(t))dt\notag\\
& + (\sqrt{\e} \lambda(\e))^{-1}\Big[\tilde{R}\Big(t,u^0(t)+\sqrt \e \lambda(\e)X^\e(t)\Big)-\tilde{R}(t,u^0(t))\Big]dt \notag\\
&=\lambda^{-1}(\e)\sigma\Big(t, u^0(t)+\sqrt\e \lambda(\e)X^{\e}(t)\Big)dW(t)+\sigma\Big(t, u^0(t)+\sqrt\e \lambda(\e)X^{\e}(t)\Big)\dot \phi^\e(t)dt,
\end{align}
with initial value $X^\e(0)=0$.

\vskip0.3cm
Now we are ready to state the second main result.
We assume that\\
\textbf{(C6)} there exist $\kappa,\ C>0$, for any $t_1,\ t_2\in [0,T]$, $u\in V$,
\begin{align}\label{C6}
|\sigma(t_1,u)-\sigma(t_2,u)|_{L_{Q}} \leq C(1+\|u\|) |t_1- t_2|^\kappa.
\end{align}
\bthm[Moderate Deviation Principle]\label{MDP}
{\rm
Under the conditions (C1-1)-(C1-3), (C2-1)-(C2-2) and (C3)-(C6),   $(u^\e-u^0)/(\sqrt\e \lambda(\e))$ obeys an LDP on $C([0,T];H)\cap L^2([0,T];V)$
with speed $\lambda^2(\e)$ and with rate function $I$ given by
\beq
I(g):=\inf_{\{h\in \HH_0;g=X^h\}}\left\{\frac12\int_0^T|\dot h(s)|_0^2ds\right\},\ \ \ \forall g\in C([0,T];H)\cap L^2([0,T];V),
\nneq
with the convention $\inf\{\emptyset\}=\infty$.
}
\nthm

We will adopt the following weak convergence method to prove the MDP.

\bthm\label{thm BD}{\rm (Budhiraja and Dupuis \refcite{Budhiraja-Dupuis}) For $\e>0$, let $\Gamma^\e$ be a measurable mapping from $C([0,T];H)$ into $\EE$.
Let $Y^\e:=\Gamma^\e(W(\cdot))$. Suppose that $\{\Gamma^\e\}_{\e>0}$ satisfies the following assumptions:
there  exists a measurable map $\Gamma^0:C([0,T];H)\rightarrow \EE$ such that
\begin{itemize}
   \item[(a)] for every $N<+\infty$ and any family $\{h^\e;\e>0\}\subset \AA_N$ satisfying that $h^\e$ converge in distribution as $S_N$-valued random elements to $h$ as $\e\rightarrow 0$,
    $\Gamma^\e\left(W(\cdot)+\frac{1}{\sqrt\e}\int_0^{\cdot}\dot h^\e(s)ds\right)$ converges in distribution to $\Gamma^0(\int_0^{\cdot}\dot h(s)ds)$ as $\e\rightarrow 0$;
   \item[(b)] for every $N<+\infty$, the set
   $
 \left\{\Gamma^0\left(\int_0^{\cdot}\dot h(s)ds\right); h\in S_N\right\}
  $
   is a compact subset of $\EE$.
 \end{itemize}
Then the family $\{Y^\e\}_{\e>0}$ satisfies a LDP in $\EE$ with the rate function $I$ given by
\beq\label{rate function}
I(g):=\inf_{\{h\in \HH_0;g=\Gamma^0(\int_0^{\cdot}\dot h(s)ds)\}}\left\{\frac12\int_0^T|\dot h(s)|_0^2ds\right\},\ g\in\EE,
\nneq
with the convention $\inf\{\emptyset\}=\infty$.
 }\nthm

\vskip0.3cm
For $h\in \HH_0$, set
$\Gamma^0\left(\int_0^\cdot \dot h(s)ds\right):=X^h$ and define $\Gamma^\e$ satisfies
$\Gamma^{\e}\left(W(\cdot)+\frac{\lambda(\e)}{\sqrt\e}\int_0^{\cdot}\dot \phi^\e(s)ds\right)=X^\e.$
To prove Theorem \ref{MDP}, we only need to verify the following two propositions according to Theorem \ref{thm BD}.
 \bprop\label{Prop weak convergence}{\rm Under the same conditions as Theorem \ref{MDP}, for every fixed  $N\in\mathbb{N}$, let $\phi^\e,\ \phi\in \mathcal{A}_N$ be such that $\phi^\e$ convergence in
distribution to $\phi$ as $\e\rightarrow0$. Then
\begin{center}
    $\Gamma^\e\left(W(\cdot)+\frac{\lambda(\e)}{\sqrt\e}\int_0^{\cdot}\dot \phi^\e(s)ds\right)$ convergence in
distribution to $\Gamma^0\left(\int_0^{\cdot}\dot \phi(s)ds\right)$
\end{center}in $C([0,T];H)\cap L^2([0,T];V)$ as $\e\rightarrow0$.}
\nprop
\bprop\label{Prop Gamm 0 compact}{\rm Under the same conditions as Theorem \ref{MDP},
  for every positive number $N<\infty$,
  the family
  $K_N:= \left\{\Gamma^0\left(\int_0^{\cdot}\dot h(s)ds\right); h\in S_N\right\}$
  is compact in $C([0,T];H)\cap L^2([0,T];V)$.}
 \nprop

We start to prove Proposition \ref{Prop weak convergence} and we need the following two lemmas. In this section, set $B(u,v,w)=(B(u,v),w)$.
\begin{lem}\label{Prop skeleton}{\rm Under the same conditions as Theorem \ref{MDP}, for any $h\in \HH_0$ and $\phi^\e\in\mathcal{A}$,  \eqref{eq skeleton} and \eqref{eq main} admit the
 unique solutions, respectively, $X^h$, $X^\e$ in $C([0,T];H)\cap L^2([0,T];V)$. Moreover, for any $N>0$ and $p\ge 2$, there exist constants $c_{T, N}$ and $\e_0>0$ such that for any $h\in S_N$, $\phi^\e\in \AA_N$,
\beq\label{eq sobolev}
\sup_{\e\in(0,\e_0]}\Big[\ee\left(\sup_{0\le t\le T}|X^\e(t)|^{2p}\right)+\ee\left(\int_0^T|X^\e(s)|^{2p-2}
\|X^\e(s)\|^{2}ds \right)\Big]
\le c_{T, N}
\nneq
and specially,
\beq\label{eq skeleton estimate}
\sup_{0\le t\le T}|X^h(t)|^{2p}+\int_0^T|X^h(s)|^{2p-2}\|X^h(s)\|^{2}ds
\le c_{T, N}.
\nneq
}
\end{lem}

\bprf\ \ Similarly as Theorem 2.4 in \refcite{CM}, the existence and uniqueness of the solution can be proved.
Here, we will prove \eqref{eq sobolev} and \eqref{eq skeleton estimate}.

Define $\tau_M:=\inf\{t: |X^{\e}(t)|^2+\int_0^t\|X^{\e}(s)\|^2ds> M\}$. By It\^o's formula and \eqref{eq main}, we have
\begin{align*}
&|X^{\e}(t\wedge\tau_M)|^2+2\int_0^{t\wedge\tau_M}\|X^{\e}(s)\|^2ds
+2\int_0^{t\wedge\tau_M}\left(B(X^\e(s),u^0(s)), X^\e(s)\right)ds \\
&=- 2(\sqrt{\e} \lambda(\e))^{-1}\int_0^{t\wedge\tau_M}\left(\tilde{R}(s,u^0(s)+\sqrt \e \lambda(\e)X^\e(t))-\tilde{R}(s,u^0(s)),  X^\e(s)\right)ds \\
&+
2\lambda^{-1}(\e)\int_0^{t\wedge\tau_M}\Big(\sigma\left(s,u^0(s)+\sqrt\e \lambda(\e)X^\e(s)\right)dW(s),X^\e(s)\Big)\\
&+2\int_0^{t\wedge\tau_M}\Big(\sigma\left(s,u^0(s)+\sqrt\e \lambda(\e)X^\e(s) \right)\dot \phi^\e(s), X^\e(s)\Big)ds\\
&+\lambda^{-2}(\e)\int_0^{t\wedge\tau_M}|\sigma\left(s,u^0(s)+\sqrt\e \lambda(\e)X^\e(s)\right)|_{{L_Q}}^2ds.
\end{align*}
Applying It\^o's formula to $f(x)=x^p$ with $x=|X^{\e}(t\wedge\tau_M)|^2$, we have $df(x)=px^{p-1}dx+\frac12 p(p-1)x^{p-2}d<x>$. Then
\begin{align}\label{00001}
&|X^{\e}(t\wedge\tau_M)|^{2p}\notag\\
=&\int_0^{t\wedge\tau_M} p|X^{\e}(s)|^{2(p-1)}d|X^{\e}(s)|^{2}
+\frac{1}{2}p(p-1)  \int_0^{t\wedge\tau_M}|X^{\e}(s)|^{2(p-2)}d<|X^{\e}|^{2}>(s) \notag\\
=&-2p\int_0^{t\wedge\tau_M}|X^{\e}(s)|^{2(p-1)}\|X^{\e}(s)\|^2ds
-2p\int_0^{t\wedge\tau_M}|X^{\e}(s)|^{2(p-1)}\left(B(X^\e(s),u^0(s)), X^\e(s)\right)ds \notag\\
& -2p(\sqrt{\e} \lambda(\e))^{-1} \int_0^{t\wedge\tau_M} |X^{\e}(s)|^{2(p-1)}\left(\tilde{R}(s,u^0(s)+\sqrt \e \lambda(\e)X^\e(t))-\tilde{R}(s,u^0(s)),  X^\e(s)\right)ds \notag\\
&+
2p\lambda^{-1}(\e)\int_0^{t\wedge\tau_M}|X^{\e}(s)|^{2(p-1)}\Big(\sigma\left(s,u^0(s)+\sqrt\e \lambda(\e)X^\e(s)\right)dW(s),X^\e(s)\Big)
\notag\\
&+2p\int_0^{t\wedge\tau_M}|X^{\e}(s)|^{2(p-1)}\Big(\sigma\left(s,u^0(s)+\sqrt\e \lambda(\e)X^\e(s) \right)\dot \phi^\e(s), X^\e(s)\Big)ds
\notag\\
&+p\lambda^{-2}(\e)\int_0^{t\wedge\tau_M}|X^{\e}(s)|^{2(p-1)}|\sigma\left(s,u^0(s)+\sqrt\e \lambda(\e)X^\e(s)\right)|_{{L_Q}}^2ds
\notag\\
&+
\frac12p(p-1) \cdot 4 \lambda^{-2}(\e)\int_0^{t\wedge\tau_M}|X^{\e}(s)|^{2(p-2)}
\Big|\Big(\sigma\left(s,u^0(s)+\sqrt\e \lambda(\e)X^\e(s)\right),X^\e(s)\Big)\Big|_{L_Q(H_0,R)}^2ds
\notag\\
:=&\sum\limits^{7}_{k=1}I_k.
\end{align}
By \eqref{hydro07}, one has
\begin{align}\label{I101}
|I_2(t)|&\le  2p\int_0^{t\wedge\tau_M}|X^{\e}(s)|^{2(p-1)} \left(\eta\|X^{\e}(s)\|^2+ C_{\eta}|X^{\e}(s)|^2\cdot\|u^0(s)\|^4_{\mathcal{H}}\right)ds.
\end{align}

From (C3), we have
\begin{align}\label{I1102}
|I_3(t)|&\le 2 p R_1
\int_0^{t\wedge\tau_M}|X^{\e}(s)|^{2p} ds.
\end{align}

By the Burkholder-Davis-Gundy inequality and (C2-1), we have
 \begin{align}\label{I203}
&\ee\left(\sup_{0\le s\le t}|I_4(s)|\right)\notag\\
=&
2p\lambda^{-1}(\e)\ee \left(\sup_{0\le s\le t} \left|\int_0^{s\wedge \tau_M }|X^{\e}(l)|^{2(p-1)}\Big(\sigma\left(l,u^0(l)+\sqrt\e \lambda(\e)X^\e(l)\right)dW(l),X^\e(l)\Big)\right|\right)\notag\\
\le& 2p\lambda^{-1}(\e) \ee\left(\int_0^{t\wedge\tau_M}
|X^{\e}(s)|^{2(p-1)\times2+2}
|\sigma\left(s,u^0(s)+\sqrt\e \lambda(\e)X^\e(s)\right)|_{L_Q}^2ds\right)^{1/2}\notag\\
\le& 2p\lambda^{-1}(\e) \ee \Big(\sup_{0\le s\le t}
|X^{\e}(s\wedge \tau_M)|^{p}  \notag\\
 &\cdot \Big[\int_0^{t\wedge \tau_M}|X^{\e}(s)|^{2p-2}\left(K_0+ 2K_1|u^0(s)|^2+2K_1\e \lambda^2(\e)|X^\e(s)|^2\right)ds\Big]^{\frac12}\Big)
\notag\\
\le& 1/4 \ee\left(\sup_{0\le s\le t}
|X^{\e}(s\wedge \tau_M)|^{2p}  \right)
 +
 c_p\lambda^{-2}(\e) \ee\left(\int_0^{t\wedge\tau_M}
|X^{\e}(s)|^{2p-2}\left(K_0+ 2K_1|u^0(s)|^2\right)ds  \right)
\notag\\
& + c_{p,K_1}\e \ee\left(\int_0^{t\wedge\tau_M}
|X^{\e}(s)|^{2p}ds\right)
\notag\\
\le& 1/4 \ee\left(\sup_{0\le s\le t}
|X^{\e}(s\wedge \tau_M)|^{2p}  \right)
 +
 c_p\lambda^{-2}(\e) \left(K_0+ 2K_1\sup_{0\le s\le T}|u^0(s)|^2\right)\notag\\
 &\cdot \ee\left(\int_0^{t\wedge\tau_M}
1+|X^{\e}(s)|^{2p} ds  \right)
+ c_{p,K_1}\e\ee\left(\int_0^{t\wedge\tau_M}
|X^{\e}(s)|^{2p}ds  \right)
\notag\\
\le&
\Big(1/4+c_{p,K_0,K_1,T}\Big(\lambda^{-2}(\e)+\e\Big)\Big)\ee\left(\sup_{0\le s\le t}
|X^{\e}(s\wedge \tau_M)|^{2p}  \right)
+
c_{p,K_0,K_1,T}\lambda^{-2}(\e).
%2p\lambda^{-1}(\e) \ee\left(\sup_{0\le s\le t\wedge \tau_M}
%|X^{\e}(s)|^{2p}  \right)
% +
% 2p\lambda^{-1}(\e) \left(K_0+ 2K_1\cdot C(1+|\xi|^2)\right)\notag\\
% &\cdot\ee\left(\int_0^{t\wedge\tau_M}
%1+|X^{\e}(s)|^{2p} ds  \right)
%+ 4pK_1\e \lambda(\e) \ee\left(\int_0^{t\wedge\tau_M}
%|X^{\e}(s)|^{2p}ds  \right).
\end{align}
Because of
\begin{align}\label{I1003}
|\sigma\left(s,u^0(s)+\sqrt\e \lambda(\e)X^\e(s) \right)|_{L_Q}
&\le
 \left(K_0+  K_1|u^0(s)+\sqrt\e \lambda(\e)X^\e(s)|^2\right)^{\frac12}
\notag\\
&\le
 \sqrt{K_0}+  \sqrt{K_1}|u^0(s)+\sqrt\e \lambda(\e)X^\e(s)|
\notag\\
&\le
 \sqrt{K_0}+  \sqrt{K_1}|u^0(s)|+\sqrt{K_1}\sqrt\e \lambda(\e)|X^\e(s)|,
\end{align}
we get
\begin{align}\label{I1004}
&|I_5(t)|\notag\\
\le &
 2p\int_0^{t\wedge\tau_M}|X^{\e}(s)|^{2(p-1)} |X^{\e}(s)|\cdot|\sigma\left(s,u^0(s)+\sqrt\e \lambda(\e)X^\e(s) \right)|_{L_Q} | \dot{\phi}^\e(s)|ds
\notag\\
\le&
 2p\int_0^{t\wedge\tau_M}|X^{\e}(s)|^{2p-1}
 \left( \sqrt{K_0}+  \sqrt{K_1}|u^0(s)|\right)| \dot{\phi}^\e(s)|ds\notag\\
 &+2p \sqrt{K_1}\sqrt\e \lambda(\e)\int_0^{t\wedge\tau_M}|X^{\e}(s)|^{2p}
 | \dot{\phi}^\e(s)|ds
\notag\\
\le&
 2p\int_0^{t\wedge\tau_M}\left(1+|X^{\e}(s)|^{2p}\right)
 \left( \sqrt{K_0}+  \sqrt{K_1}|u^0(s)|\right)| \dot{\phi}^\e(s)|ds\notag\\
 &+2p \sqrt{K_1}\sqrt\e \lambda(\e)\int_0^{t\wedge\tau_M}|X^{\e}(s)|^{2p}
 | \dot{\phi}^\e(s)|ds
\notag\\
\le& c_{p,K_0,K_1,T,N}
 %p\int_0^{t\wedge\tau_M}\left( \sqrt{K_0}+  \sqrt{K_1}|u^0(s)|\right)^2ds+p\int_0^{t\wedge\tau_M}  | \dot{\phi}^\e(s)|^2ds
 + 2p\int_0^{t\wedge\tau_M}|X^{\e}(s)|^{2p} \left( \sqrt{K_0}+  \sqrt{K_1}|u^0(s)|\right)| \dot{\phi}^\e(s)|ds\notag\\
 &
 +c_{p,K_1,N,T}\sqrt\e \lambda(\e)\left(\sup_{0\le s\le t}
|X^{\e}(s\wedge \tau_M)|^{2p}  \right).
%+2p \sqrt{K_1}\sqrt\e \lambda(\e)\int_0^{t\wedge\tau_M}|X^{\e}(s)|^{2p}
% | \phi^\e(s)|ds.
\end{align}

\begin{align}\label{I11005}
&|I_6(t)|+|I_7(t)| \notag\\
\le &\left(p\lambda^{-2}(\e)+
2p(p-1) \lambda^{-2}(\e)\right)
\int_0^{t\wedge\tau_M}|X^{\e}(s)|^{2(p-1)}|\sigma\left(s,u^0(s)+\sqrt\e \lambda(\e)X^\e(s)\right)|_{{L_Q}}^2ds
 \notag\\
 \le &
 \left(p\lambda^{-2}(\e)+
2p(p-1) \lambda^{-2}(\e)\right)
\int_0^{t\wedge\tau_M}|X^{\e}(s)|^{2(p-1)}\left( K_0 +2K_1|u^0(s)|^2+2K_1\e \lambda^2(\e)|X^\e(s)|^2\right)ds
 \notag\\
 \le &
 c_{p,K_0,K_1,T}(\lambda^{-2}(\e)+\e)\left(\sup_{0\le s\le t}|X^{\e}(s\wedge \tau_M)|^{2p}  \right)
 +
 c_{p,K_0,K_1,T}\lambda^{-2}(\e).
 %
%
% \left(p\lambda^{-2}(\e)+
%2p(p-1) \lambda^{-2}(\e)\right)
%\int_0^{t\wedge\tau_M}\left(1+|X^{\e}(s)|^{2p}\right)\left( K_0 +2K_1|u^0(s)|^2\right)ds
%\notag\\
% &
%+\left(p\lambda^{-2}(\e)+
%2p(p-1) \lambda^{-2}(\e)\right)
%\int_0^{t\wedge\tau_M}2K_1\e \lambda^2(\e)|X^\e(s)|^{2p}ds.
\end{align}

Taking the supremum up to time $t$ in \eqref{00001}, and then taking the expectation, Gronwall's inequality implies
\begin{align}\label{3 I}
&\Big(3/4-c(\lambda^{-2}(\e)+\e+\sqrt\e\lambda(\e)\Big)
\ee\left(\sup_{0\le s\le t}|X^{\e}(s\wedge \tau_M)|^{2p}\right)\notag\\
&+\left(2p- 2p\eta\right)\ee\left(\int_0^{t\wedge\tau_M}
|X^{\e}(s)|^{2(p-1)}\|X^{\e}(s)\|^2ds\right)\notag\\
\le& C\exp\Big[{\int_0^T\Big(2pC_\eta\|u^0(s)\|^4_\mathcal{H}+2pR_1+2p(\sqrt K_0+\sqrt K_1|u^0(s)|)|\dot{\phi}^\e(s)|\Big)ds}\Big]\notag\\
\le & C_N,
\end{align}
for any $\phi^\e\in \AA_N$, here we choose $\e$ and $\eta$ small enough.
Letting $M\rightarrow \infty$, we get \eqref{eq sobolev}.

The proof of \eqref{eq skeleton estimate} is very similar to that of  the deterministic case \eqref{eq sobolev}, we omit it here.
The proof of this lemma is complete.
\nprf

\vskip 0.3cm
%%%%%%%%%%%%%%%%%%%%%%%%%%%%%%%%%%%%%%%%%%%%%%%%%%%%%%%%%%%%%%%%%%%%%%%%%%%%%%%%
For every integer $n$, let $\psi_n:\ [0,T]\rightarrow[0,T]$ denote
a measurable map such that for every $s\in [0,T]$, $s\le \psi_n(s)
\le (s+c2^{-n})\wedge T$ for some positive constant $c$. Let
$G_M^\e(t)=\{ \omega:\  (\underset{0\le s\le t}\sup |X^\e (s)|^2)
\vee  ( \int_0^t \|X^\e (s)\|^2 ds)  \le M \}$. We show the following lemma.
\begin{lem}\label{lemma01}{\rm There exist $\e_0,\ N_0>0$ such that
for any $\e\in(0,\e_0]$ and $ \phi^\e \in \AA_N$,
 \begin{align}\label{lemma4.2}
\ee \left( \mathbf{1}_{G_M^\e(T)} \int_0^T |
X^\e\big(\psi_n (s)\big) -X^\e (s)|^2 ds   \right) \le CM^3N 2^{-n/2},
\end{align}
for any $M>N_0$, where $C$ is independent on $\e$, $M$ and $N$.
}\end{lem}
\bprf\ \
Letting $\sup_{0\le s\le T} |u^0 (s)|^2+
 \int_0^T \|u^0 (s) \|^2 ds  \le N_0$. $It\hat{o}$'s formula yields
\begin{align}\label{00005}
&\ee \left( \mathbf{1}_{G_M^\e(T)} \int_0^T |X^\e\big(\psi_n (s)\big)-X^\e (s)  |^2 ds \right)  \notag\\
=& \ee \left( \mathbf{1}_{G_M^\e(T)} \int_0^T 2\int_s^{\psi_n (s)} \left( X^\e (r)
- X^\e(s), d X^\e (r) \right)  dr ds \right)\notag\\
&+ \ee \left( \mathbf{1}_{G_M^\e(T)} \int_0^T \int_s^{\psi_n (s)}\lambda^{-2}(\e)
| \sigma\left( r,\sqrt\e \lambda(\e) X^\e (r)+u^0(r)  \right) |_{L_Q}^2dr ds \right)\notag\\
=& 2 \ee \left( \mathbf{1}_{G_M^\e(T)} \int_0^T \int_s^{\psi_n (s)} \left( X^\e (r)- X^\e(s), -A X^\e (r) \right)drds  \right)\notag\\
&-2 \ee \left( \mathbf{1}_{G_M^\e(T)} \int_0^T \int_s^{\psi_n (s)} \Big( X^\e (r)- X^\e(s),B \left(X^\e (r),\sqrt\e \lambda(\e) X^\e (r)+u^0(r)  \right)  \Big)drds  \right)\notag\\
&-2 \ee \left( \mathbf{1}_{G_M^\e(T)} \int_0^T \int_s^{\psi_n (s)} \Big( X^\e (r)- X^\e(s),B \left(u^0(r), X^\e (r)  \right)  \Big)drds  \right)\notag\\
&-2 \ee  \mathbf{1}_{G_M^\e(T)} \int_0^T \int_s^{\psi_n (s)} \frac{1}{\sqrt\e \lambda(\e)}\notag\\
&\cdot\left( X^\e (r)- X^\e(s),\tilde R \left(r,\sqrt\e \lambda(\e) X^\e (r)+u^0(r)  \right)- \tilde R \left(r,u^0(r)  \right) \right)drds  \notag\\
&+2 \ee \left( \mathbf{1}_{G_M^\e(T)} \int_0^T \int_s^{\psi_n (s)} \lambda^{-1}(\e)\Big( X^\e (r)- X^\e(s),\sigma \left(r,\sqrt\e \lambda(\e) X^\e (r)+u^0(r)  \right) dW(r)\Big)ds  \right)\notag\\
&+2 \ee \left( \mathbf{1}_{G_M^\e(T)} \int_0^T \int_s^{\psi_n (s)} \Big( X^\e (r)- X^\e(s),\sigma \left(r,\sqrt\e \lambda(\e) X^\e (r)+u^0(r)  \right)\dot{\phi}^\e(r)\Big) drds  \right)\notag\\
&+ \ee \left( \mathbf{1}_{G_M^\e(T)} \int_0^T \int_s^{\psi_n (s)} \lambda^{-2}(\e)
| \sigma\left( r,\sqrt\e \lambda(\e) X^\e (r)+u^0(r)  \right) |_{L_Q}^2dr ds \right)\notag\\
=:&\sum\limits^{7}_{k=1}I_k.
\end{align}
Now we estimate $I_k$, respectively.
\begin{align*}
I_1 =& 2 \ee \left( \mathbf{1}_{G_M^\e(T)} \int_0^T \int_s^{\psi_n (s)} - \|X^\e (r)   \| ^2 + \|X^\e (s)   \| \|X^\e (r)   \|drds  \right)\\
\le &  1/2\ee \left( \mathbf{1}_{G_M^\e(T)} \int_0^T  \|X^\e (s)   \| ^2 ds \int_s^{\psi_n (s)} dr  \right)\\
\le &C\cdot M \cdot 2^{-(n+1)}.
\end{align*}
For $I_2$,
\begin{align*}
I_2 \le &
2 \ee \left( \mathbf{1}_{G_M^\e(T)} \int_0^T \int_s^{\psi_n (s)} |B \left(X^\e (r), u^0(r), X^\e (r)  \right) | drds  \right)\\
&+ 2 \ee \left( \mathbf{1}_{G_M^\e(T)} \int_0^T \int_s^{\psi_n (s)} |B \left(X^\e (r),\sqrt\e \lambda(\e)X^\e (r), X^\e (s)  \right) | drds  \right)\\
&+ 2 \ee \left( \mathbf{1}_{G_M^\e(T)} \int_0^T \int_s^{\psi_n (s)} |B \left(X^\e (r), u^0(r), X^\e (s)  \right) | drds  \right)\\
=:&J_1+J_3+J_3,\ {\rm where}
\end{align*}
\begin{align*}
J_1 \le &
2 \ee \left( \mathbf{1}_{G_M^\e(T)} \int_0^T \int_s^{\psi_n (s)}
\eta \| X^\e (r)  \|^2 +C_{\eta}| X^\e (r)  |^2 \| u^0(r)  \|_{\mathcal{H}}^4  drds  \right)\\
\le& 2 \eta\ee \left( \mathbf{1}_{G_M^\e(T)} \int_0^T\| X^\e (r)  \|^2 dr \int_{(r-c2^{-n})\vee0}^{r}
ds  \right)\\
&+ 2 C_{\eta} a_0^2 \ee \left( \mathbf{1}_{G_M^\e(T)} \int_0^T \int_s^{\psi_n (s)}
| X^\e (r)  |^2\cdot| u^0(r)  |^2\| u^0(r)  \|^2   drds  \right)\\
\le& 2\eta C2^{-n}M
+ 2 C_{\eta}\cdot a_0^2 \cdot M^2 \ee \left( \mathbf{1}_{G_M^\e(T)} \int_0^T \int_s^{\psi_n (s)}
 \| u^0(r)  \|^2  drds  \right)\\
\le& 2\eta C2^{-n}M
+ 2 C_{\eta}\cdot a_0^2 \cdot M^2 \ee \left( \mathbf{1}_{G_M^\e(T)} \int_0^T dr \int_{(r-c2^{-n})\vee0}^{r}  \| u^0(r)  \|^2 ds  \right)\\
\le& 2\eta CM2^{-n}
+ 2 C_{\eta}\cdot a_0^2 \cdot M^3 2^{-n};
\end{align*}
\begin{align*}
J_2 \le &
2\sqrt\e \lambda(\e) \ee \left( \mathbf{1}_{G_M^\e(T)} \int_0^T \int_s^{\psi_n (s)}
\eta \| X^\e (r)  \|^2 +C_{\eta}| X^\e (r)  |^2 \| X^\e (s)  \|_{\mathcal{H}}^4  drds  \right)\\
\le& 2 \sqrt\e \lambda(\e) \eta CM2^{-n}
+ 2 \sqrt\e \lambda(\e) C_{\eta} \cdot M
 \ee \left( \mathbf{1}_{G_M^\e(T)} \int_0^T \| X^\e (s)  \|_{\mathcal{H}}^4  ds \right)\cdot 2^{-n}\\
\le& 2 \sqrt\e \lambda(\e) \eta CM2^{-n}
+ 2 \sqrt\e \lambda(\e) C_{\eta}a_0^2 \cdot M
 \ee \left( \mathbf{1}_{G_M^\e(T)}\sup_{0\le s\le T}| X^\e (s)  |^2 \int_0^T   \| X^\e (s)  \|^2  ds \right)\cdot 2^{-n}\\
\le& 2 \sqrt\e \lambda(\e) \eta CM2^{-n}
+ 2\sqrt\e \lambda(\e) C_{\eta}\cdot a_0^2 \cdot M^3 2^{-n};
\end{align*}
by (\ref{hydro04}), we have
\begin{align*}
J_3 \le &
2 \ee \left( \mathbf{1}_{G_M^\e(T)} \int_0^T \int_s^{\psi_n (s)}
\eta \| X^\e (s)  \|^2 +C_{\eta}\| X^\e (r)  \|_{\mathcal{H}}^2 \| u^0 (r)  \|_{\mathcal{H}}^2  drds  \right)\\
\le& 2  CM2^{-n}
+ 2 a_0^2 C_{\eta}
 \ee \left( \mathbf{1}_{G_M^\e(T)} \int_0^T \int_s^{\psi_n (s)}
 | X^\e (r)  | \cdot   \| X^\e (r)  \| \cdot
 | u^0 (r)  | \cdot   \| u^0 (r)  \| drds \right)\\
\le& 2  CM2^{-n}
+ 2  C_{\eta}a_0^2 M
 \ee \left( \mathbf{1}_{G_M^\e(T)} \int_0^T \int_s^{\psi_n (s)}
    \| X^\e (r)  \|
     \| u^0 (r)  \| drds \right)\\
\le& 2  CM2^{-n}
+  C_{\eta}a_0^2 M
 \ee \left( \mathbf{1}_{G_M^\e(T)} \int_0^T \int_s^{\psi_n (s)}
    \| X^\e (r)  \|^2+
     \| u^0 (r)  \|^2 drds \right)\\
\le& 2  CM2^{-n}
+ C_{\eta}a_0^2 M^22^{-n},
\end{align*}
here integral transformation as in $J_1$ has been applied to obtain the last inequality.

Similarly as $J_3$, we have
\begin{align*}
I_3 \le &
2 \ee \left( \mathbf{1}_{G_M^\e(T)} \int_0^T \int_s^{\psi_n (s)} |\left (B (u^0(r), X^\e (r)), X^\e (s)   \right) | drds  \right)\\
\le &
2\ee \left( \mathbf{1}_{G_M^\e(T)} \int_0^T \int_s^{\psi_n (s)}
\eta \| X^\e (s)  \|^2 +C_{\eta}\| X^\e (r)  \|_{\mathcal{H}}^2 \| u^0 (r)  \|_{\mathcal{H}}^2  drds  \right)\\
\le& CM2^{-n}
+ C_{\eta}a_0^2 M^22^{-n}.
\end{align*}
It is easy to see
\begin{align*}
I_4 \le &
2 \ee \left( \mathbf{1}_{G_M^\e(T)} \int_0^T \int_s^{\psi_n (s)}
|X^\e (r)- X^\e (s)    |\cdot R_1 |X^\e (r)    |drds  \right)\\
\le & C_TM2^{-n}.
\end{align*}
Because of ${G_M^\e(T)}\subset {G_M^\e(r)} $, we have
\begin{align*}
I_5 \le &
 2\lambda^{-1}(\e)\int_0^T  \ee \left( \int_s^{\psi_n (s)} \mathbf{1}_{G_M^\e(r)}| X^\e (r)- X^\e(s)|^2\cdot |\sigma \left(r,\sqrt\e \lambda(\e) X^\e (r)+u^0(r)  \right)|_{L_Q}^2dr  \right)^{\frac12} ds\notag\\
 \le &
 2\lambda^{-1}(\e)\int_0^T  \ee \left( \int_s^{\psi_n (s)} \mathbf{1}_{G_M^\e(r)} | X^\e (r)- X^\e(s)|^2\cdot \left(K_0+K_1|\sqrt\e \lambda(\e) X^\e (r)+u^0(r)|^2  \right)dr  \right)^{\frac12} ds\notag\\
\le &  \lambda^{-1}(\e)  CM2^{-\frac n2} ;
\end{align*}
\begin{align*}
I_6 \le &
2\ee \left( \mathbf{1}_{G_M^\e(T)}\   \int_0^T  \int_s^{\psi_n (s)} | X^\e (r)- X^\e(s)|\cdot |\sigma \left(r,\sqrt\e \lambda(\e) X^\e (r)+u^0(r)  \right)|_{L_Q} \cdot   |\dot{\phi}^\e(r)| drds \right)\notag\\
 \le & CM
  \ee \left( \mathbf{1}_{G_M^\e(T)}\   \int_0^T  \int_s^{\psi_n (s)}   |\dot{\phi}^\e(r)| drds \right)\notag\\
  \le & CM
  \ee \left( \mathbf{1}_{G_M^\e(T)}\   \int_0^T
   \left( \int_s^{\psi_n (s)}   |\dot{\phi}^\e(r)|^2 dr
    \cdot \int_s^{\psi_n (s)}   dr
     \right)^\frac12
    ds \right)\notag\\
\le & C_TM \cdot N^\frac12  \cdot 2^{-\frac n2}.
\end{align*}
Recall (\ref{h}), we have for all sufficient small $\e$,
\begin{align*}
I_7 \le &
\ee \left( \mathbf{1}_{G_M^\e(T)}\   \int_0^T  \int_s^{\psi_n (s)}
\frac{1}{\lambda^2(\e)} \left(K_0+K_1|\sqrt\e \lambda(\e) X^\e (r)+u^0(r)|^2  \right)drds
\right)\notag\\
\le & C_TM  \cdot 2^{-n}.
\end{align*}
Putting all the estimation of $I_k$ into \eqref{00005}, we obtain
\eqref{lemma4.2} which completes the proof.
\nprf
\vskip 0.3cm

{\bf Proof of Propsition\ref{Prop weak convergence}:}
From Skorokhod representation theorem, there exists $(\tilde \phi_\e, \tilde \phi, \tilde W^\e)$ such that $(\tilde \phi_\e, \tilde W^\e)$ and $( \phi^\e,  W)$ have the same distribution, the distribution of $\tilde \phi$ coincides with that of $\phi$ and $\tilde \phi_\e$ converges $\tilde \phi$ a.s. in the weak topology of $S_M$. To lighten notations, we will write $(\tilde \phi_\e, \tilde \phi, \tilde W^\e)=(\phi^\e, \phi, W)$.
\vskip 0.2cm

Recall (\ref{eq skeleton}) and (\ref{eq main}). Let $Z_\e=X^\e-X^\phi$, here $X^\phi$ is the solution of (\ref{eq skeleton}) replaced $h$ by $\phi$, then $Z_\e(0)=0$ and
\begin{align*}
&dZ_{\e}(t)+AZ_{\e}(t)dt+B\Big(X^{\e}(t), \sqrt \e \lambda(\e)X^\e(t)\Big)dt+B( Z_\e(t),u^0(t)) dt+B(u^0(t), Z_\e(t)) dt\notag\\
& + \left((\sqrt{\e} \lambda(\e))^{-1}[\tilde{R}(t,u^0(t)+\sqrt \e \lambda(\e)X^\e(t))-\tilde{R}(t,u^0(t))]-\tilde R'(t,u^0(t))X^\phi(t) \right)dt\notag\\
=&\lambda^{-1}(\e)\sigma\Big(t, u^0(t)+\sqrt\e \lambda(\e)X^{\e}(t)\Big)dW(t)+\sigma\Big(t, u^0(t)+\sqrt\e \lambda(\e)X^{\e}(t)\Big)\dot \phi^\e(t)dt\notag\\
&
- \sigma\Big(t, u^0(t)\Big)\dot \phi(t)dt .
\end{align*}
Hence
\begin{align*}
&|Z_{\e}(t)|^2+2\int_0^{t}\|Z_{\e}(s)\|^2ds\\
=&2\int_0^{t} \sqrt\e \lambda(\e) B\left(X^\e(s), X^\e(s),X^\phi(s)\right)ds
-2\int_0^{t}B\left(Z_\e(s),u^0(s), Z_\e(s)\right)ds \\
&- 2\int_0^{t}\left((\sqrt{\e} \lambda(\e))^{-1}[\tilde{R}(s,u^0(s)+\sqrt \e \lambda(\e)X^\e(s))-\tilde{R}(s,u^0(s))]-\tilde R'(s,u^0(s))X^\phi(s),  Z_\e(s)\right)ds \\
&+
2\lambda^{-1}(\e)\int_0^{t}\Big(\sigma\left(s,u^0(s)+\sqrt\e \lambda(\e)X^\e(s)\right)dW(s),Z_\e(s)\Big)\\
&+2\int_0^{t}\Big(\sigma\left(s,u^0(s)+\sqrt\e \lambda(\e)X^\e(s) \right)\cdot \dot{\phi}^\e(s)-\sigma\Big(s, u^0(s)\Big)\cdot \dot{\phi}(s), Z_\e(s)\Big)ds\\
&+\lambda^{-2}(\e)\int_0^{t}|\sigma\left(s,u^0(s)+\sqrt\e \lambda(\e)X^\e(s)\right)|_{{L_Q}}^2ds\\
\le&2 \sqrt\e \lambda(\e)\int_0^{t} \eta \|X^\e(s)\|^2 +C_\eta |X^\e(s)|^2 \|X^\phi(s)\|_{\mathcal{H}}^4ds
+2\int_0^{t} \eta  \|Z_\e(s) \|^2+ C_\eta |Z^\e(s)|^2 \|u^0(s)\|_{\mathcal{H}}^4 ds \\
&+2\int_0^{t}\Big|(\sqrt{\e} \lambda(\e))^{-1}\Big[\tilde{R}\Big(s,u^0(s)+\sqrt \e \lambda(\e)X^\e(s)\Big)-\tilde{R}(s,u^0(s)) -\tilde R'(s,u^0(s))
\sqrt{\e} \lambda(\e)X^\e(s) \Big]\Big| \notag\\
&\cdot |Z_\e(s)|ds \\
&+ 2\int_0^{t}|\tilde R'(s,u^0(s))| \cdot |Z_\e(s)|^2ds
+
2\lambda^{-1}(\e)\int_0^{t}\Big(\sigma\left(s,u^0(s)+\sqrt\e \lambda(\e)X^\e(s)\right)dW(s),Z_\e(s)\Big)\\
&+2\int_0^{t}\Big( \sigma\Big(s, u^0(s))(\dot{\phi}^\e(s)-\dot{\phi}(s)),Z_\e(s)\Big)ds
+2\int_0^{t} \sqrt{L_1} \sqrt\e \lambda(\e) |X^\e(s)|\cdot |\dot{\phi}^\e(s)| \cdot |Z_\e(s)| ds\\
&+\lambda^{-2}(\e)\int_0^{t}K_0+K_1|u^0(s)+\sqrt\e \lambda(\e)X^\e(s)|^2ds\\
:=& \sum\limits^{8}_{k=1}I_k(t).
\end{align*}
Since,
 $$I_7 \le \e \lambda^2(\e) \int_0^{t} {L_1}|X^\e(s)|^2 |\dot{\phi}^\e(s)|^2ds  + \int_0^{t} |Z_\e(s)|^2 ds:=I'_7  + \int_0^{t} |Z_\e(s)|^2 ds,$$
 and
\begin{align*}
I_3 =& 2\int_0^{t}\Big|\int_0^1\Big[\tilde{R}'\Big(s,u^0(s)+x\sqrt \e \lambda(\e)X^\e(s)\Big)-\tilde{R}'(s,u^0(s))\Big]X^\e(s)dx\Big||Z_\e(s)|ds\\
\le & C\sqrt\e \lambda(\e)\int_0^t|X^\e(s)|^2|Z_\e(s)|ds\\
\le& C\e \lambda^2(\e)\int_0^t|X^\e(s)|^4ds+\int_0^{t} |Z_\e(s)|^2ds \\
:=& I'_3 + \int_0^{t}|Z_\e(s)|^2ds.
\end{align*}
Then,
\begin{align*}
&|Z_{\e}(t)|^2+2(1-\eta)\int_0^{t}\|Z_{\e}(s)\|^2ds\\
\le& I_1+I'_2+I'_3+I_5+I_6+I'_7+I_8 +\int_0^{t}\left(C_\eta \|u^0(s)\|_{\mathcal{H}}^4+2+2 |\tilde R'(s,u^0(s))| \right)|Z_\e(s)|^2ds.
\end{align*}
%here
%$$
%I'_1=C_\eta \sqrt\e \lambda(\e)\int_0^{t} |X^\e(s)|^2 \|X^\phi(s)\|_{\mathcal{H}}^4ds.
%$$
Gronwall's inequality implies
\begin{align}\label{0050}
&|Z_{\e}(t)|^2+2(1-\eta)\int_0^{t}\|Z_{\e}(s)\|^2ds \notag\\
\le& (I_1+I'_3+I_5+I_6+I'_7+I_8)\exp\Big(\int_0^{T}C_\eta \|u^0(s)\|_{\mathcal{H}}^4+2+2(\tilde{R}'_1+\tilde{R}'_2|u^0(s)|)ds\Big).
\end{align}

By Lemma \ref{Lem 1} and (\ref{eq skeleton estimate}), we may take $M$ big enough such that $\underset{h\in S_N}\sup \Big(\underset{0\le s\le T}\sup |X^h (s)|^2+\int_0^T \|X^h (s)\|^2ds\Big)\bigvee\Big(\underset{0\le s\le T}\sup |u^0 (s)|^2+\int_0^T \|u^0 (s)\|^2ds\Big)\le M$. Let
$G_{M,\e}(t)=\{\underset{0\le s\le t}\sup |X^\e (s)|^2)\le M\}
\cap  \{  \int_0^t \|X^\e (s) \|^2 ds)  \le M \}$. We split the proof into two steps.
\vskip 0.2cm

Step 1: $\pp\Big(\Big(G_{M,\e}(T)\Big)^c\Big)\le C(1+|\xi|^4)M^{-1}$ which is obtained by Chebysev inequality and (\ref{eq sobolev}).

Step 2: $\ee \Big( \mathbf{1}_{G_{M,\e}(T)}( \underset{0\le s\le T}\sup |Z_{\e} (s)|^2+\int_0^T \|Z_{\e} (s)\|^2ds  ) \Big)\rightarrow 0$ as $\e \rightarrow0$.

Now we show the Step 2. When $\omega \in G_{M,\e}(T)$, we have\\
$$I_1 \le C\sqrt{\e} \lambda(\e) \cdot M^3,\ \ \ I'_3 \le C_T\epsilon\lambda^2(\e)M^2.$$
For $I_5$,
\begin{align*}
\ee\left(\mathbf{1}_{G_{M,\e}(T)}\underset{0\le t\le T}\sup | I_5(t)| \right)
&\le 2\lambda^{-1}(\e) \ee\left(\int_0^{T} \mathbf{1}_{G_{M,\e}(t)}
|\sigma\left(t,u^0(t)+\sqrt\e \lambda(\e)X^\e(t)\right)|_{L_Q}^2|Z_\e(t)|^2dt \right)^{\frac12}\\
&\le C_T \lambda^{-1}(\e)  M.
\end{align*}
Applying Lemma \ref{lemma01}, similar to (3.23) in \refcite{CM},
$$\ee\left(\mathbf{1}_{G_{M,\e}(T)}\underset{0\le t\le T}\sup  I_6(t,\e) \right)\rightarrow0,\ \ \  \e \rightarrow0.$$
And
$$I'_7 \le C\e \lambda^2(\e)\cdot N\cdot M, \ \ \ I_8 \le C_T\lambda^{-2}(\e)(1+M).$$

Thus, $\underset{\e \rightarrow0}\lim \ee \big( \mathbf{1}_{G_{M,\e}(T)}( \underset{0\le s\le T}\sup |Z_{\e} (s)|^2+\int_0^T \|Z_{\e} (s)\|^2ds  ) \big)=0$. So Step 2 is obtained.

Finally, since for any $\delta>0$,
\begin{align*}
& \mathbb{P}\Big(\underset{0\le s\le T}\sup |Z_{\e} (s)|^2+\int_0^T \|Z_{\e} (s)\|^2ds>\delta\Big)\\
\le &
\pp\Big(\Big(G_{M,\e}(T)\Big)^c\Big)
+1/\delta \ee \Big( \mathbf{1}_{G_{M,\e}(T)}( \underset{0\le s\le T}\sup |Z_{\e} (s)|^2+\int_0^T \|Z_{\e} (s)\|^2ds  ) \Big),
\end{align*}
Step 1 and Step 2 implies that $\lim_{\e\rightarrow0}\mathbb{P}\Big(\underset{0\le s\le T}\sup |Z_{\e} (s)|^2+\int_0^T \|Z_{\e} (s)\|^2ds>\delta\Big)=0$, which concludes the proof of the proposition.

\vskip 0.3cm

{\bf Proof of Propsition\ref{Prop Gamm 0 compact}:} The proof is similar as that of Proposition \ref{Prop weak convergence} and easier. The proof will be omitted.

\appendix
%
%\section{Appendices}
%{Appendices should be used only when absolutely necessary. They
%should come\hfilneg}
%\eject
%\noindent
%before the Acknowledgment. If there is more than one
%appendix, number them alphabetically. Number displayed equations
%in the way, e.g.~(\ref{that}), (A.2), etc.
%\begin{equation}
%f(j\delta, i\delta) \cong \frac{\pi}{M} \sum^M_{n=1}
%Q_{\theta_n} (j \cos \theta_n + i \sin \theta_n)\,. \label{that}
%\end{equation}
%
\section*{Acknowledgment}
This work was supported by National Natural Science Foundation of China (NSFC) (No. 11431014, No. 11401557), and the Fundamental Research Funds for the Central Universities (No. 0010000048).

%\section*{Citations}
%They are to be cited in the text in superscript
%after the punctuation marks e.g.~word,\refcite{2} and word:\refcite{2}. If it is mentioned in the text as part of a sentence, it should be of normal size, e.g.~see Ref.~\refcite{3}.

\end{document}